\documentclass[a4paper]{article}
\usepackage{cancel}
\usepackage{indentfirst}
\usepackage{amsmath}
\usepackage{amssymb}
\usepackage{amsthm}
\usepackage{upgreek}
\usepackage{graphicx}
\usepackage{varioref}
\usepackage{hyperref}
\usepackage{mathtools}
\usepackage{mathrsfs}
\usepackage{changepage}
\usepackage{nicefrac}
\usepackage[normalem]{ulem}
\usepackage{tensor}
\usepackage{stmaryrd}
\usepackage{empheq}
\usepackage{xstring}
\usepackage{cleveref}
\usepackage{trimclip}
\usepackage{centernot}
\usepackage{xcolor}
\usepackage{enumerate}
\usepackage{wasysym}
\usepackage[scr=boondox]{mathalpha}
\usepackage{tikz}
\usepackage{tikz-cd}
\usetikzlibrary{decorations.pathmorphing}
\usetikzlibrary{calc}
\usetikzlibrary{babel}

\tikzset{curve/.style={settings={#1},to path={(\tikztostart)
.. controls ($(\tikztostart)!\pv{pos}!(\tikztotarget)!\pv{height}!270:(\tikztotarget)$)
and ($(\tikztostart)!1-\pv{pos}!(\tikztotarget)!\pv{height}!270:(\tikztotarget)$)
.. (\tikztotarget)\tikztonodes}},
    settings/.code={\tikzset{quiver/.cd,#1}
    \def\pv##1{\pgfkeysvalueof{/tikz/quiver/##1}}},
    quiver/.cd,pos/.initial=0.35,height/.initial=0}

\tikzset{tail reversed/.code={\pgfsetarrowsstart{tikzcd to}}}
\tikzset{2tail/.code={\pgfsetarrowsstart{Implies[reversed]}}}
\tikzset{2tail reversed/.code={\pgfsetarrowsstart{Implies}}}
\tikzset{no body/.style={/tikz/dash pattern=on 0 off 1mm}}

\DeclareMathAlphabet{\mathmybb}{U}{bbold}{m}{n}
\newcommand{\1}{\mathmybb{1}}

\graphicspath{ {../pics/} }

\theoremstyle{plain}
\newtheorem{lemma_env}{\textup{\bfseries{Lemma}}}[section]
\newtheorem{proposal_env}{\textup{\bfseries{Proposal}}}[section]
\newtheorem{problem_env}{\textup{\bfseries{Problem}}}[section]
\newtheorem{corollary_env}{\textup{\bfseries{Corollary}}}[section]
\newtheorem{fact_env}{\textup{\bfseries{Fact}}}[section]
\newtheorem{statement_env}{\textup{\bfseries{Statement}}}[section]
\newtheorem{exercise_env}{\textup{\bfseries{Exercise}}}[section]

\theoremstyle{remark}
\newtheorem*{note_env}{Note}
\newtheorem*{supplement_env}{Supplement}
\newtheorem*{example_env}{Example}
\newtheorem*{examples_env}{Examples}
\newtheorem*{properties_env}{Properties}
\newtheorem*{precaution_env}{Precaution}
\newtheorem*{counterexample_env}{Counterexample}
\newtheorem*{counterexamples_env}{Counterexamples}
\newtheorem*{ifact_env}{Interesting Fact}
\newtheorem*{question_env}{Question}

\theoremstyle{definition}
\newtheorem{theorem_env}{\textup{\bfseries{Theorem}}}[section]
\newtheorem{definition_env}{\textup{\bfseries{Definition}}}[section]

\crefname{lemma_env}{Lemma}{Lemmata}
\crefname{proposal_env}{Proposal}{Proposals}
\crefname{problem_env}{Problem}{Problems}
\crefname{corollary_env}{Corollary}{Corollaries}
\crefname{fact_env}{Fact}{Facts}
\crefname{statement_env}{Statement}{Statements}
\crefname{exercise_env}{Exercise}{Exercises}
\crefname{theorem_env}{Theorem}{Theorems}
\crefname{definition_env}{Definition}{Definitions}
\crefname{section}{Section}{Sections}
\crefname{subsection}{Subsection}{Subsections}
\crefname{figure}{Figure}{Figures}

\newcommand{\theorem}[2][]{\begin{theorem_env}[#1]#2\end{theorem_env}}
\newcommand{\lemma}[2][]{\begin{lemma_env}[#1]#2\end{lemma_env}}
\newcommand{\corollary}[2][]{\begin{corollary_env}[#1]#2\end{corollary_env}}

\newcommand{\proposal}[2][]{\begin{proposal_env}[#1]#2\end{proposal_env}}

\newcommand{\all}[1]{\begin{cases}#1\end{cases}}

\newcommand{\any}[1]{\left[\begin{aligned}#1\end{aligned}\right.}

\newcommand{\prove}[2][Proof]{\begin{proof}[#1]#2\end{proof}}
\newcommand{\provebullets}[2][Proof]{\prove[#1]{\bullets{#2\qedhere}}}

\renewcommand{\gather}[1]{\begin{gather*}#1\end{gather*}}
\renewcommand{\multline}[1]{\begin{multline*}#1\end{multline*}}

\newcommand{\bullets}[1]{\begin{itemize}#1\end{itemize}}
\newcommand{\numbers}[1]{\begin{enumerate}#1\end{enumerate}}
\newcommand{\defset}[2]{\left\{#1\mid#2\right\}}

\newcommand{\ind}[3]{\tensor[#1]{#2}{#3}}

\newcommand{\bigboxplus}{
    \mathop{
        \vphantom{\bigoplus}
        \mathchoice
        {\vcenter{\hbox{\resizebox{\widthof{$\displaystyle\bigoplus$}}{!}{$\boxplus$}}}}
        {\vcenter{\hbox{\resizebox{\widthof{$\bigoplus$}}{!}{$\boxplus$}}}}
        {\vcenter{\hbox{\resizebox{\widthof{$\scriptstyle\oplus$}}{!}{$\boxplus$}}}}
        {\vcenter{\hbox{\resizebox{\widthof{$\scriptscriptstyle\oplus$}}{!}{$\boxplus$}}}}
    }\displaylimits
}

\newcommand{\proddots}{\cdot\ldots\cdot}

\newcommand{\Z}{\mathbb{Z}}

\renewcommand{\O}{\mathrm{O}}
\newcommand{\M}{\mathrm{M}}

\renewcommand{\ge}{\geqslant}
\renewcommand{\le}{\leqslant}

\newcommand{\map}{\rightarrow}
\newcommand{\Map}{\longrightarrow}

\newcommand{\sm}{\setminus}

\newcommand{\normeq}{\trianglelefteqslant}

\newcommand{\vect}[1]{\begin{pmatrix}#1\end{pmatrix}}

\newcommand{\smallvect}[1]{\begin{psmallmatrix}#1\end{psmallmatrix}}

\newcommand{\arr}[2]{\begin{array}{#1}#2\end{array}}

\newcommand{\angles}[1]{\left\langle#1\right\rangle}

\renewcommand{\tilde}{\widetilde}

\newcommand{\fm}{\mathfrak{m}}

\renewcommand{\phi}{\varphi}

\newcommand{\cat}[1]{\StrChar{#1}{1}[\firstchar]\StrGobbleLeft{#1}{1}[\remainingchars]\mathscr{\firstchar}\IfStrEq{\remainingchars}{}{}{\!\cat{\remainingchars}}}

\DeclareMathOperator{\id}{id}

\DeclareMathOperator{\Ker}{Ker}

\DeclareMathOperator{\Image}{Im}

\DeclareMathOperator{\End}{End}

\DeclareMathOperator{\Cent}{Cent}

\DeclareMathOperator{\GL}{GL}
\DeclareMathOperator{\EO}{EO}

\DeclareMathOperator{\Trans}{Trans}

\newcounter{lection}
\date{}

\title{Normal Structure of Isotropic Odd Orthogonal Groups}
\author{
    Leonid Danilevich \\
    Saint Petersburg University, \\
    7/9 Universitetskaya nab., \\
    St. Petersburg, 199034 Russia
}

\begin{document}

    \maketitle
    \begin{abstract}
        Let $(M, q)$ be a quadratic projective module of an odd rank over an commutative ring, where the form $q$ is semiregular, with global Witt index of at least $2$, and with $\mathrm{rk}(M) \ge 7$.
        We prove standard commutator formulae and classify $\EO$-normal subgroups of $\O(M, q)$ without assumption of $2$ being invertible.
    \end{abstract}

    \section{Introduction}

    The normal structure of the general linear groups $\GL(n, K)$ over a commutative ring $K$ for $n \ge 3$ is described by the Wilson–Golubchik theorem~\cite{wilson-golubchik, golubchik}. Namely, if $H \le \GL(n, K)$ is normalized by the elementary subgroup $\mathrm{E}(n, K)$, then there exists a unique ideal $I \normeq K$ such that
    $\mathrm{E}(n, K, I) \le H \le \mathrm{C}(n, K, I)$, where on the left is the relative elementary subgroup, and on the right is the full congruence subgroup.

    Subsequently, similar results were obtained by L. Vaserstein for arbitrary Chevalley groups under the assumption of invertibility of the structure constants~\cite{vasserstein} and later by E. Abe in the general case~\cite{abe} for all root systems of rank $\ge 3$.
    For the remaining root systems $\mathrm{B}_2$ and $\mathrm{G}_2$, complete answers in the split case were obtained by D. Costa and G. Keller~\cite{costa, keller}.

    In the rank $1$ case, such a classification cannot be obtained in general, since the classification of normal subgroups of $\mathrm{SL}(2, \Z)$ essentially includes the classification of all finite groups generated by two elements of orders $2$ and $3$.
    Nevertheless, under additional arithmetic assumptions, the answer was found by Costa and Keller~\cite{costakeller}.

    If the structure constants are invertible, then the levels in the normal structure are parameterized by ideals of the ring, as in the case of general linear groups.
    If the root system is multiply laced and is not $\mathrm{B}_2$, then they are parameterized by ``admissible pairs'' in the sense of Abe, and for $\mathrm{B}_2$ (and $\mathrm{A}_1$ in the context of the article by Costa and Keller) so-called radices appear.

    V. Petrov and A. Stavrova in~\cite{petrov-stavrova} began to study the groups of points of isotropic reductive group schemes over arbitrary commutative rings.
    In this article, they gave the definition of an elementary subgroup and proved its normality.
    Later, Stavrova and A. Stepanov~\cite{stavrova-stepanov} found the normal structure under the assumption of invertibility of the structure constants.

    Finally, the normal structure is also known for the class of unitary groups, which includes classical isotropic reductive groups. The most general result in this direction was obtained by A. Bak and R. Preusser in~\cite{buck-preusser} for odd unitary groups of rank $\ge 3$ without the assumption of invertibility of $2$.

    This work considers the simplest of the remaining cases of isotropic reductive groups.
    Namely, the normal structure is obtained for orthogonal groups $\O(M, q)$, where $M$ is a projective module of odd rank with a semiregular quadratic form $q$, the global Witt index of $q$ is at least $2$ and the rank of $M$ is at least $7$ (the case of rank $5$ is covered by the article of Costa and Keller).
    As expected, the levels are parameterized by admissible pairs in the sense of Abe.

    The proof of the standard classification is based on a variant of the work of Preusser~\cite{preusser}.
    In addition to the main result (Theorem~\ref{main-theorem}), we will also obtain
    \bullets{
        \item normality of the elementary subgroup $\EO(M, q) \normeq \O(M, q)$ for an arbitrary quadratic form $q$ of Witt index $2$ with nontrivial anisotropic kernel~(\cref{eo-normality});
        \item standard commutator formulae and an analogue of perfection for relative elementary subgroups in the case of a semiregular quadratic form~(\cref{standard-commutator,eo-perfection});
        \item explicit computation of the transporter involved in the normal structure (Proposition~\ref{transporters}).
    }

    \section{Basic Definitions}
    For group elements $x, y \in G$, their left-normed commutator is denoted $[x, y] = xyx^{-1}y^{-1}$, and conjugation of $x$ by $y$ on the left is $\ind{^y}{x}{} = yxy^{-1}$.

    All rings are assumed to be commutative with unity.
    The group of invertible elements of a ring $K$ is denoted by $K^\times$.

    Let $M$ be an arbitrary $K$-module.
    Its dual module of $K$-linear functionals is denoted by $M^\lor$.

    Everywhere it is assumed that $M$ is equipped with a fixed quadratic form — a map $q: M \to K$, homogeneous of degree $2$ ($\forall u \in M$, $x \in K$: $q(ux) = q(u)x^2$), such that its polarization $b_q: M \times M \to K$, $u, v \mapsto q(u + v) - q(u) - q(v)$ is bilinear.
    The pair $(M, q)$ will be called a quadratic module.

    For an element $u \in M$, the $K$-linear functional $u^* \in M^\lor$ is defined, $u^* = b_q(u, \_)$.
    Accordingly, $M^* = \defset{u^*}{u \in M} \le M^\lor$. 

    Everywhere it is assumed that $(M, q)$ has Witt index at least $2$, i.e., $M$ can be represented as $M = e_{-2}K \oplus e_{-1}K \oplus M_0 \oplus e_1 K \oplus e_2 K$, and $q(x_{-2}, x_{-1}, u_0, x_1, x_2) = q_0(u_0) + x_{-1}x_1 + x_{-2}x_2$, where $q_{0}$ is a quadratic form on $M_0$.
    This decomposition of $M$ into a direct sum is assumed to be fixed and will not change.

    For an arbitrary $K$-linear endomorphism $\phi: M \to M$, its matrix is denoted by $[\phi]$: $[\phi]_{ij} = \pi_i\phi \iota_j$ for $i,j \in \{0, \pm 1,\pm 2\}$.
    Here $\iota_0: M_0 \to M$, $\pi_0: M \to M_0$ are the canonical embedding and projection, similarly for $\iota_i: K \to M$ and $\pi_i: M \to K$ for $i \ne 0$.
    The rows and columns are numbered in order $(-2, -1, 0, 1, 2)$.
    Zero elements are omitted, and irrelevant ones are denoted by $*$.

    The orthogonal group $\mathrm{O}(M, q)$ is the subgroup of the group $\GL(M)$ of $K$-linear automorphisms of $M$, consisting of automorphisms that preserve the quadratic form.

    If $M$ is a projective module, one can define $\mathrm{SO}(M, q) = \O(M, q) \cap \mathrm{SL}(M)$.
    Equivalently, $\mathrm{SO}(M, q)$ consists of those automorphisms $\phi \in \O(M, q)$ such that in the localization at any maximal ideal $\fm \normeq K$, the corresponding automorphism $\phi_{\fm}$ of the free module $M_{\fm}$ has determinant one.

    In this case, the decomposition holds:
    \[\O(M, q) = \mathrm{SO}(M, q) \times \mu_2(K)\id,\]
    where $\mu_2(K) = \defset{\omega \in K}{\omega^2 = 1}$.

    \subsection{Semiregular Quadratic Form}
    Assume that $M_0$ is a projective module of odd rank $n \ge 1$.
    Denote $n = 2\nu + 1$.

    Let $\xi_i, \xi_{ij}$ be formal variables, with $\xi_{ij} = \xi_{ji}$.
    Construct the matrix $\Xi$ of size $n\times n$, filled with them according to the formula $\Xi_{ij} = \all{2\xi_i,&i=j \\ \xi_{ij}, &i \ne j}$.

    As is known~\cite[IV.3.1.2]{knus}, $\det(\Xi) \in \Z[\xi_i, \xi_{ij}]$ is divisible by $2$; define $P_n = \frac{1}{2}\det(\Xi) \in \Z[\xi_i, \xi_{ij}]$.
    The form $q_0$ is said to be \emph{semiregular} if the ideal $(P_n(q_0(u_i), b_{q_0}(u_i, u_j)) \mid u_1, \ldots, u_n \in M_0) = (1) \normeq K$.

    We will need the following proposition about the structure of a quadratic module with a semiregular form over a local ring~\cite[IV.3.1.7]{knus}:
    \proposal{\label{orthogonal-decomposition}
    Let $K$ be a local ring.
    For every free semiregular quadratic $K$-module $(M_0, q_0)$ of rank $2\nu + 1$, there exists an orthogonal decomposition
        \[(M_0, q_0) = \bigboxplus\limits_{i = 1}^{\nu}(M_i', q'_i) \boxplus (M'_0, q'_{0}),\]
        where $(M'_{i}, q'_{i})$ for $1 \le i \le \nu$ is a free quadratic module of rank $2$ with a nondegenerate quadratic form, and $(M'_{0}, q'_{0})$ is a free quadratic module of rank $1$ with a semiregular quadratic form.
        The module $M'_{i}$ has a basis $\{f_{-i}, f_i\}$ such that $b_{q'_{i}}(f_{-i}, f_i) = 1$ and $M'_0$ has a basis $\{f_0\}$ such that $q'_{0}(f_0) \in K^\times$.
    }
    The following two lemmas are also easy to see:
    \lemma{\label{unit-q}
    If $(M_0, q_0)$ is a projective module of odd rank $n$ with a semiregular quadratic form, then the ideal $(q_0(u)\mid u \in M_0) = (1) \normeq K$.
    \prove{
        The semidiscriminant $\frac{1}{2}\det(\Xi)$ is a polynomial with integer coefficients in $q_0(u)$ for $u \in M_0$.
    }
    }
    \lemma{\label{unit-bq}
    If $(M_0, q_0)$ is a projective module of odd rank $n \ge 3$ with a semiregular quadratic form, then the ideal $(b_{q_0}(u, v)\mid u, v \in M_0) = (1) \normeq K$.
    \prove{
        It suffices to show that the semidiscriminant $\frac{1}{2}\det(\Xi) \in (2\xi_i, \xi_{ij}) \normeq \Z[\xi_i, \xi_{ij}]$.

        Write the semidiscriminant $\frac{1}{2}\det(\Xi)$ as a linear combination of monomials in $\xi_i, \xi_{ij}$.
        All monomials except $\xi_1\proddots \xi_n$ lie in the ideal generated by $\xi_{ij}$, and the latter enters with coefficient $2^{n}$, and since $n \ge 3$, $\frac{1}{2}(2^n\xi_1\proddots \xi_n)$ still lies in the ideal $(2\xi_1)$.
    }
    }

    \subsection{Elementary Subgroup}
    \lemma{
        The following endomorphisms of $M$ lie in $\O(M, q)$:
        \bullets{
            \item $T_{ij}(y) = 1 + e_{ij}y - e_{-j,-i}y$, where $y \in K$, $i, j \in \{\pm 1, \pm 2\}$, and $i \ne \pm j$;
            \[[T_{1,2}(y)] = \vect{1&-y&&&\\&1&&&\\&&1&&\\&&&1&y\\&&&&1};\]
            \item $T_{i}(u) = 1 + e_{0,i}u - e_{-i,0}u^* - e_{-i,i}q_0(u)$, where $u \in M_0$, $i \in \{\pm 1, \pm 2\}$;
            \[[T_{1}(u)] = \vect{1&&&&\\&1&-u^*&-q_0(u)&\\&&1&u&\\&&&1&\\&&&&1};\]
            \item $d_i(r) = 1 + e_{ii}(r - 1) + e_{-i,-i}(r^{-1} - 1)$, where $r \in K^\times$, $i \in \{\pm 1, \pm 2\}$;
            \[[d_{1}(r)] = \vect{1&&&&\\&r^{-1}&&&\\&&1&&\\&&&r&\\&&&&1};\]
            \item $d(\phi_0) = 1 + e_{0,0}(\phi_0 - 1)$, where $\phi_0 \in O(M_0, q_0)$.
            \[[d(\phi_0)] = \vect{1&&&&\\&1&&&\\&&\phi_0&&\\&&&1&\\&&&&1}.\]
        }
    }
    The elementary subgroup $\EO(M, q) = \angles{T_{ij}(y), T_i(u)\mid y \in K, u \in M_0}$ is a subgroup of $\O(M, q)$.
    Its generators $T_{ij}(y)$ and $T_i(u)$ will be called \emph{elementary transvections}.

    Subgroups of $\O(M, q)$ normalized by $\EO(M, q)$ will be called $\EO$-normal.

    The diagonal subgroup $\mathrm{D}(M, q) = \angles{d_i(r), d(\phi_0)\mid r \in K^\times, \phi_0 \in \O(M_0, q_0)}$ is another subgroup of $\O(M, q)$.

    Let $L$ be a $K$-algebra.
    One can define the quadratic $L$-module $(M_L, q_L)$:
    as an $L$-module, it is $M \otimes_K L$, and the quadratic form $q_L$ is defined naturally:
    \[q_L\left(\sum\limits_{i = 1}^{N}m_i \otimes \ell_i\right) = \sum\limits_{i=1}^{N}q(m_i)\ell_i^2 + \sum\limits_{1 \le i < j \le N}b_{q}(m_i, m_j)\ell_i \ell_j.\]
    It is easy to check that this definition is well-formed.

    For the quadratic module $(M_L, q_L)$, one can take the orthogonal group $\O(M_L, q_L)$, and this is a functor from the category of $K$-algebras to the category of groups: for every arrow $L \to L'$ between $K$-algebras, $M_{L'} = M_L \otimes_L L'$, and for $\phi \in \O(M_L, q_L)$, $\phi \otimes_L L' \in \O(M_{L'}, q_{L'})$ is well-defined.

    Let $I \normeq K$ be an ideal.
    Then one can define the relative orthogonal subgroup $\O(M, q, I)$ as the kernel of $\O(M, q) \map \O(M/MI, q_{K/I})$.
    Also define the relative elementary subgroup $\mathrm{FO}(M, q, I) = \angles{T_{ij}(y), T_i(u)\mid y \in I, u \in M_0 I}$ for some ideal $I \normeq K$.
    This subgroup is almost never normal in $\O(M, q)$.
    The correct definition of the relative elementary subgroup is given after Lemma~\ref{some-lemma}.
    Clearly, $\mathrm{FO}(M, q, I) \le \O(M, q, I)$.

    \subsection{Elementary Relations}
    \lemma{
        For $y, z \in K$, $u, v \in M_0$, $r \in K^\times$, $\phi_0 \in O(M_0, q_0)$, $i, j \in \{\pm 1, \pm 2\}$ with $i \ne \pm j$:
        \begin{align*}
            T_{ij}(y)\cdot  T_{ij}(z) &= T_{ij}(y + z),\\
            T_{i}(u)\cdot  T_{i}(v) &= T_{i}(u + v),\\
            [T_{i}(u), T_{j}(v)] &= T_{-j,i}(b_{q_0}(u, v)),\\
            [T_i(u), T_{ij}(y)] &= T_j(u y) \cdot T_{-j,i}(q_0(u) y), \\
            [T_i(u), T_{ji}(y)] &= 1, \\
            [T_{ij}(y), T_{i,-j}(z)] &= 1;
        \end{align*}
        Other reasonable commutation formulae between elementary transvections follow from the identities $T_{ij}(y) = T_{-j,-i}(-y)$;
        \begin{align*}
            \ind{^{d(\phi_0)}}{T_{ij}(y)}{} &= T_{ij}(y),\\
            \ind{^{d(\phi_0)}}{T_{i}(u)}{} &= T_{i}(\phi_0 u),\\
            \ind{^{d_i(r)}}{T_{ij}(y)}{} &= T_{ij}(y r), \\
            \ind{^{d_j(r)}}{T_{ij}(y)}{} &= T_{ij}(y r^{-1}), \\
            \ind{^{d_i(r)}}{T_{i}(u)}{} &= T_{i}(u r^{-1}), \\
            \ind{^{d_j(r)}}{T_{i}(u)}{} &= T_{i}(u).
        \end{align*}
        Since $d_{i}(r) = d_{-i}(r^{-1})$, other commutation formulae between diagonal matrices and elementary transvections follow trivially.
    }

    \section{Gauss Decomposition}
    \lemma{\label{moving-invertibles}
    \bullets{Let $(K, \fm)$ be a local ring, and $\phi \in \O(M, q)$.
    \item[a.] Then $\exists \psi \in \EO(M, q)$ such that $[\phi\cdot\psi]_{-1,-1} \in K^\times$.
    \item[b.] Suppose $\phi$ fixes $e_{-2}$ and $e_2$.
    Then $\psi$ can be chosen so that it also fixes $e_{-2}$ and $e_2$.
    }
    \prove{
        Since $\phi \in \GL(M)$, $\Image(\phi) = M$.
        In particular, $\Image(\pi_{-1} \circ \phi) = K$, and the ideal\\ $([\phi]_{-1,-2},[\phi]_{-1,-1}, [\phi]_{-1,1},[\phi]_{-1,2}) + [\phi]_{-1,0}(M_0) \normeq K$ is the unit ideal.
        \bullets{
            \item If $[\phi]_{-1,-1} \in K^\times$, then $\psi \coloneqq \id$ works.
            \item Else if $\exists j \in \{\pm 2\}$: $[\phi]_{-1,j} \in K^\times$, then $\psi \coloneqq T_{-1,j}(1)$ works (this is case (a), since in case (b) it is guaranteed that $[\phi]_{-1,j} = 0$).
            \item Else if $\exists u \in M_0$ such that $[\phi]_{-1,0}\cdot u - [\phi]_{-1,1}\cdot q(u)$ is invertible, then $\psi\coloneqq T_{-1}(u)$ works.
            \item Else we obtain that for any $u \in M_0$: $[\phi]_{-1,0}\cdot u - [\phi]_{-1,1}\cdot q(u) \in \fm$. However, since the ideal generated by the terms is the unit ideal, $\exists w \in M_0$: both $[\phi]_{-1,0}\cdot w, [\phi]_{-1,1}\cdot q(w) \in K^\times$.
            Substituting $u \coloneqq w \cdot x$, we get the identity $\forall x \in K$: \[[\phi]_{-1,0}\cdot w \cdot x - [\phi]_{-1,1}\cdot q(w)\cdot x^2 \equiv 0\pmod{\fm}.\]
            This is a nonzero polynomial of degree $2$ in $x$, vanishing everywhere on $K/\fm$, hence $K/\fm$ is the field of two elements.
            At the same time, $[\phi]_{-1,1}\cdot q(w) \in K^\times$, so $q(w) \equiv 1\pmod{\fm}$.
            Let's compute $\ind{^{T_{-1}(w)}}{T_{1}(w)}{}$, depicting in the matrix only the rows and columns with indices $-1, 0, 1$ (other elements coincide with the corresponding elements of the identity matrix):
            \multline{\left[\ind{^{T_{-1}(w)}}{T_{1}(w)}{}\right] =  \vect{1 &  &  \\ w& 1 & &\\-q_0(w) &-w^*& 1}\cdot\vect{1 & -w^* & -q(w) \\ & 1 & w & \\ && 1}\cdot\vect{1 &  &  \\ -w& 1 &  & \\ -q_0(w)&w^*& 1} =\\= \vect{(1 + q(w))^2 & -(1 + q(w))w^* & -q(w) \\ w (q(w)^2 + q(w)) & 1 - q(w) w w^* & w(1 - q(w)) \\ -q(w)^3 & (q(w)^2 - q(w))w^* & (1 - q(w))^2} \equiv  \vect{&&1\\ &*&\\1&&}\pmod{\fm}.}
            Thus $\psi \coloneqq \ind{^{T_{-1}(w)}}{T_{1}(w)}{}$ works.\qedhere
        }
    }
    }
    \lemma{\label{bottom-row-is-zero}
    Let $\phi \in \O(M, q)$ send $e_{-2}$ to $e_{-2} + e_2 y$ for some $y \in K$.
    Then $\phi$ has the form
        \[[\phi] = \vect{1 & * & * & * & *\\ &*&*&*&*\\ &*&*&*&*\\ &*&*&*&*\\ &&&&1}.\]
        \provebullets{
            \item Since $\phi$ preserves the quadratic form, $y = q(\phi(e_{-2})) = q(e_{-2}) = 0$.
            \item Further, $\phi(e_{\pm 1})$ is orthogonal to $\phi(e_{-2})$, so $[\phi]_{2,\pm1} = 0$.
            \item Moreover, $\forall u \in M_0$: $\phi(u)$ is orthogonal to $\phi(e_{-2})$, so $\forall u \in M_0: [\phi]_{2,0}\cdot u = 0$.
            Hence, $[\phi]_{2,0} = 0 \in M^\lor$.
            \item Finally, $[\phi]_{2,2} = b_q(\phi(e_{-2}), \phi(e_2)) = b_q(e_{-2}, e_2) = 1$.\qedhere
        }
    }

    \proposal{
        Let $(K, \fm)$ be a local ring.
        Then $\mathrm{O}(M, q) = \EO(M, q) \cdot \mathrm{D}(M, q)$.

        \prove{
            From the commutation formulae, $[\EO(M, q), \mathrm{D}(M, q)] \le \EO(M, q)$.
            Hence, it suffices to reduce $\phi \in \O(M, q)$ to $\id$ by multiplying on the left and right by elements of $\EO(M, q)$ and $\mathrm{D}(M, q)$.
            \numbers{
                \item By Lemma~\ref{moving-invertibles} for position $(-2, -2)$ instead of $(-1, -1)$, we can assume, multiplying $\phi$ on the right by an element of $\EO(M, q)$, that $[\phi]_{-2,-2} \in K^\times$.
                \item Additionally, we can assume, multiplying $\phi$ (say, on the right, it doesn't matter) by $d_{-2}\left([\phi]_{-2,-2}^{-1}\right)$, that $[\phi]_{-2,-2} = 1$:
                \[[\phi] = \vect{1 & * & * & * & *\\ *&*&*&*&*\\ *&*&*&*&*\\ *&*&*&*&*\\ *&*&*&*&*}.\]
                \item Additionally, we can assume, multiplying $\phi$ on the left by elements from $\EO(M, q) \cdot \mathrm{D}(M, q)$, that $[\phi]_{i,-2} = 0$ for $i \in \{0, \pm 1\}$:
                \[[\phi] = \vect{1 & * & * & * & *\\ &*&*&*&*\\ &*&*&*&*\\ &*&*&*&*\\ *&*&*&*&*}.\]
                For this, it suffices to replace $\phi$ by $T_{-1,-2}(-[\phi]_{-1,-2})\cdot T_{0,-2}(-[\phi]_{0,-2})\cdot T_{1,-2}(-[\phi]_{1,-2})\cdot \phi$.
                \item Now by Lemma~\ref{bottom-row-is-zero}: $[\phi]_{2,j} = \delta_{2,j}$ for $j \in \{0,\pm 1,\pm 2\}$:
                \[[\phi] = \vect{1 & * & * & * & *\\ &*&*&*&*\\ &*&*&*&*\\ &*&*&*&*\\ &&&&1}.\]
                \item Similarly to step 3, but for the right column, we can additionally assume, multiplying $\phi$ on the left by elements from $\EO(M, q) \cdot \mathrm{D}(M, q)$, that $[\phi]_{i,2} = 0$ for $i \in \{0, \pm 1\}$:
                \[[\phi] = \vect{1 & * & * & * & *\\ &*&*&*&\\ &*&*&*&\\ &*&*&*&\\ &&&&1}.\]
                \item Now by the symmetric version of Lemma~\ref{bottom-row-is-zero}, $[\phi]_{-2,j} = \delta_{-2,j}$ for $j \in \{0, \pm 1, \pm 2\}$:
                \[[\phi] = \vect{1 &  &  &  & \\ &*&*&*&\\ &*&*&*&\\ &*&*&*&\\ &&&&1}\]
                \item Similarly to steps 1.--6., but for the smaller square, we can assume that nonzero values occur only on the main diagonal:
                \[[\phi] = \vect{1 &  &  &  & \\ &1&&&\\ &&*&&\\ &&&1&\\ &&&&1}\]
                This means $\phi \in \mathrm{D}(M, q)$.\qedhere
            }
        }
    }

    \section{Normality of $\EO(M, q) \normeq \O(M, q)$}

    \subsection{Continuity of Conjugation}
%
    \proposal{\label{conjugation-is-continuous}
    Let $S \subset K$ be a multiplicative system, $\phi \in \EO(M_{S^{-1}K}, q_{S^{-1}K}) \cdot\\ \cdot \mathrm{D}(M_{S^{-1}K}, q_{S^{-1}K})$.
    Also assume that the ideal $(q_0(w)\mid w \in M_0) \normeq K$ is the unit ideal, and the module $M_0$ is finitely generated.

    Then $\forall s \in S$: $\exists s' \in S$: $\ind{^\phi}{\Image(\mathrm{FO}(M, q, K s'))}{} \le \Image(\mathrm{FO}(M, q, K s))$, where the image is taken inside $\EO(M_{S^{-1}K}, q_{S^{-1}K})$.
    \prove{
        Without loss of generality, $\phi = d_i(r)$, or $\phi = d(\phi_0)$, or $\phi = T_i(u)$, or $\phi = T_{ij}(y)$.
        \bullets{
            \item Let $\phi = d_i(\nicefrac{r_1}{t})$, and $\phi^{-1} = d_i(\nicefrac{r_2}{t})$, where $r_1, r_2 \in K, t \in S$.
            From the commutation formulae, it is clear that $s' = st$ works.
            \item Let $\phi = d(\phi_0)$.
            Denote by $\tilde{M_0}$ the image of $M_0$ inside $S^{-1}M_0$, this is a $K$-submodule.
            Let $t \in S$ be such that $\phi_0(\tilde{M_0} t) \le \tilde{M_0}$ — exists by $\tilde{M_0}$ being finitely generated.
            Clearly, $\phi_0(\tilde{M_0} st) = \phi_0(\tilde{M_0} t)s \le \tilde{M_0} s$. From the commutation formulae, it is clear that $s' = st$ works.
            \item Let $\phi = T_{-j}(\nicefrac{u}{t})$, where $u \in M_0$, $t \in S$. From the commutation formulae, it is clear that (for $i \ne \pm j$) \[\ind{^\phi}{\Image(\angles{T_{kl}(K s t^2), T_{i}(M_0 s t), T_{-j}(M_0 s t))}}{} \le \Image(\mathrm{FO}(M, q, Ks)).\]
            It remains to note that $T_{j}(v) = [T_{i}(\nicefrac{v}{st^2}), T_{ij}(st^2)]T_{-j,i}(-\nicefrac{q_0(v)}{st^2})$, so $s' = s^2 t^3$ works. 
            \item Let $\phi = T_{i,-j}(\nicefrac{y}{t})$, where $y \in K$, $t \in S$. From the commutation formulae, it is clear that \[\ind{^\phi}{\Image(\angles{T_{ij}(K s), T_{i,-j}(K s), T_{-i,-j}(K s), T_{k}(M_0 s t)})}{} \le \Image(\mathrm{FO}(M, q, Ks)).\]
            Note that $T_{-j,i}(q_0(u)z) = T_{j}(-uz)\cdot[T_i(u), T_{ij}(z)]$.
            Now use the fact that the ideal $(q_0(w)\mid w \in M_0) = (1) \normeq K$, i.e., $\exists w_1, \ldots, w_N \in M_0$, $x_1, \ldots, x_N \in K$: $q_0(w_1)x_1 + \ldots + q_0(w_N)x_N = 1$. Write
            \[T_{-j,i}(z) = \prod\limits_{k = 1}^{N}T_{-j,i}(q_0(w_k)x_k z) = \prod\limits_{k = 1}^{N}T_{j}(-\nicefrac{w_k x_k z}{st})[T_i(w_k st), T_{ij}(\nicefrac{x_k z}{s^2 t^2})].\]
            It is clear that $s' = s^3 t^2$ works.\qedhere
        }
    }
    }

    \subsection{Injectivity of Localization}
    \lemma{\label{localization-is-injective}
    Let $R$ be a Noetherian ring, $M$ a finitely generated $R$-module with a quadratic form $q$, $S \subset R$ a multiplicative system.
    It is claimed that $\exists s \in S$: the homomorphism $\mathrm{O}(M, q, Rs) \to \O(S^{-1}M, q_{S^{-1}R})$ is injective.
    \prove{
        Clearly, $\mathrm{O}(M, q, Rs) = \mathrm{O}(M, q) \cap \GL(M, Rs)$ (and $\GL(M, I) = \Ker(\GL(M) \to \GL(M/I))$).
        Check that $\GL(M, Rs) \to \GL(S^{-1}M)$ is injective for $s \gg 1$.
        As is known, $\Ker(M \to S^{-1}M) = \bigcup\limits_{t \in S}(0_M : t)$, where $0_M \le M$ is the zero submodule.

        By Noetherianness, $\exists s \in S$: $\Ker(M \to S^{-1}M) = (0_M : s)$.
        Check that $(0_M : s) \cap Ms = 0_M$: indeed, for $u \in (0_M : s) \cap Ms$: $\exists v \in M: u = vs$, so $vs^2 = 0$, and by the definition of $s$, we get $0 = vs = u$.

        Now let $\phi \in \Ker(\GL(M, Rs) \to \GL(S^{-1}M))$.
        This means that $\forall u \in M$: $\phi(u) = u \in S^{-1}M$, and simultaneously $\forall u \in M$: $\phi(u) = u\pmod{Ms}$.
        According to the above, $\phi(u) = u \in M$.
    }
    }

    \subsection{Proof of the Normality}
    \theorem{\label{eo-normality}
    Let the ideal $(q_0(w)\mid w \in M_0) \normeq K$ be the unit ideal, and $M_0$ a finitely generated $K$-module.
    Then $\EO(M, q) \normeq \O(M, q)$.
    \prove{
        Let $\phi \in \O(M, q)$. Fix $y \in K, u \in M_0$.
        Check that $\ind{^\phi}{T_{ij}(y)}{}, \ind{^\phi}{T_{i}(u)}{} \in \EO(M, q)$ for all $i, j \in \{\pm 1, \pm 2\}$, $i \ne\pm j$.

        Let $M_0 = \angles{u_1, \ldots, u_m}_K$.
        Let $V = \{e_{-2}, e_{-1}, u_1, \ldots, u_m, e_{1}, e_2\}$.
        In particular, $M = \angles{V}_K$.

        Without loss of generality, $K$ is Noetherian. Namely, generate the $\Z$-algebra $R \le K$ by the following elements of $K$:
        \bullets{
            \item For each $v \in V$: $\phi(v), \phi^{-1}(v) \in M$.
            Decompose them in the generating system $V$, and add the decomposition coefficients from $K$ to $R$.
            \item Since $(q_0(w)\mid w \in M_0) = (1) \normeq K$, there exist $w_1, \ldots, w_N \in M_0$, $x_1, \ldots, x_N \in K$ such that $q_0(w_1)x_1 + \ldots + q_0(w_N)x_N = 1$.
            Each $w_i$ decomposes in the generators $u_1, \ldots, u_m$ with coefficients from $K$, add these coefficients to $R$. Also add $x_1, \ldots, x_N$ to $R$.
            \item Add to $R$ the coefficients of the decomposition of $u$ in the generators $u_1, \ldots, u_m$, and also add $y$.
            \item Add to $R$ the numbers $q_0(u_1), \ldots, q_0(u_m)$, and the numbers $b_{q_0}(u_i, u_j)$ for $1 \le i < j \le m$.
        }
        Now $R$ is Noetherian, and $M_R = \angles{e_{-2}, e_{-1}, u_1, \ldots, u_m, e_1, e_2}_R \le M$ is a finitely generated $R$-submodule.
        The restriction $q_R: M_R \times M_R \to R$ is still a quadratic form, actually taking values in $R$ by the construction of $R$.
        The restriction $\phi_R$ on $M_R$ is still an automorphism of the module $M_R$, and, since $q_R$ is the restriction of $q$, $\phi$ preserves the quadratic form, i.e., $\phi_R \in \O(M_R, q_R)$.

        As soon as it is proved that $\ind{^{\phi_R}}{T_{ij}(y)}{}, \ind{^{\phi_R}}{T_{i}(u)}{} \in \EO(M_R, q_R)$, the analogous statement for $K$ will immediately be true: say, if $\ind{^{\phi_R}}{T_{i}(u)}{} = T_{k_1}(u_1)\cdot T_{i_1,j_1}(y_1)\proddots T_{k_n}(u_n)\cdot T_{i_n,j_n}(y_n)$, then the analogous decomposition holds in $M$: both sides of the equality are $K$-linear endomorphisms of $M$, whose actions coincide on the generators. From now on, assume $M = M_R, q = q_R$.

        Consider the ideal $I = \defset{z \in R}{\ind{^{\phi_R}}{T_{ij}(Rz)}{}, \ind{^{\phi_R}}{T_{i}(M_0 z)}{} \le \EO(M, q)}$.

        Let $\fm \normeq R$ be maximal. Consider $S \coloneqq R \sm \fm$. By Lemma~\ref{localization-is-injective} on the injectivity of localization, $\exists s \in S$: $\O(M, q, R s) \hookrightarrow \O(S^{-1}M, S^{-1}R)$.
        By Proposition~\ref{conjugation-is-continuous} on the continuity of conjugation \[\exists s' \in S:\quad \ind{^{\phi_R}}{T_{ij}(R s')}{}, \ind{^{\phi_R}}{T_{i}(M_0 s')}{} \in \Image(\mathrm{FO}(M, q, R s)) \le \O(S^{-1}M, q_{S^{-1}R}).\]
        By the choice of $s$, $\ind{^{\phi_R}}{T_{ij}(R s')}{}, \ind{^{\phi_R}}{T_{i}(M_0 s')}{} \in \mathrm{FO}(M, q, R s) \le \EO(M, q)$.

        Thus, $s' \in I$, i.e., $\forall \fm \normeq R$: $I \nsubseteq \fm$. In other words, $I = (1)$.
    }
    }

    \section{Relativization}
    In this section, we define relative subgroups with respect to admissible pairs and generalize the normality theorem for the elementary subgroup to relative subgroups, using doubles of rings and modules.
    The idea of this approach goes back to Stein~\cite{stein}.

    \subsection{Doubling of the Ring and Module Along an Admissible Pair}
    Let $I \normeq K$ be an ideal, and $N_0 \le M_0$ a submodule.
    Assume that the pair $(N_0, I)$ is \emph{admissible}, namely, the conditions hold
    \[\all{M_0 I \le N_0\\ b_{q_0}(M_0, N_0) + q_0(N_0) \subset I}\label{admissible}\tag{$B$}\]
    One can define the double of the ring $K$ along the ideal $I$ and the double of the module $M_0$ along the submodule $N_0$, as the pullbacks of the corresponding diagrams:
    \[\begin{tikzcd}
    {K\times_I K}
          & K && {M_0\times_{N_0}M_0} & {M_0} \\
          K & {K/I} && {M_0} & {M_0/N_0}
          \arrow["{\pi_1}"', from=1-1, to=1-2]
          \arrow["{\pi_2}", from=1-1, to=2-1]
          \arrow["\lrcorner"{anchor=center, pos=0.125}, draw=none, from=1-1, to=2-2]
          \arrow["\pi"', from=1-2, to=2-2]
          \arrow["{\pi_1}"', from=1-4, to=1-5]
          \arrow["{\pi_2}", from=1-4, to=2-4]
          \arrow["\lrcorner"{anchor=center, pos=0.125}, draw=none, from=1-4, to=2-5]
          \arrow["\pi"', from=1-5, to=2-5]
          \arrow["\pi", from=2-1, to=2-2]
          \arrow["\pi", from=2-4, to=2-5]
    \end{tikzcd}\]
    It is easy to see that \begin{align*}
                               K\times_I K =& \defset{(x_1, x_2) \in K \times K}{x_1 \equiv x_2 \pmod{I}}, \\M_0 \times_{N_0} M_0 =& \defset{(u_1, u_2) \in M_0 \times M_0}{u_1 \equiv u_2 \pmod{N_0}}.
    \end{align*}
    Define the diagonal embeddings $\Delta: K \to K\times_I K$ and $\Delta: M_0\to M_0\times_{N_0}M_0$.

    Since $M_0 I \le N_0$, the set $M_0 \times_{N_0} M_0$ has the structure of a $K\times_{I} K$-module with componentwise action.

    Since $b_{q_0}(M_0, N_0) + q_0(N_0) \subset I$, the quadratic form $q_0$ is well-defined on $M_0 \times_{N_0} M_0$ by the formula $q_0(m_1, m_2) = (q_0(m_1), q_0(m_2))$.

    Denote $N \coloneqq e_{-2}I \oplus e_{-1}I \oplus N_0 \oplus e_1 I \oplus e_2 I$.
    Clearly, there is an isomorphism \[M \times_N M \cong e_{-2}(K\times_I K) \oplus e_{-1}(K\times_I K) \oplus (M_0 \times_{N_0} M_0 )\oplus e_1 (K\times_I K) \oplus e_2 (K\times_I K).\]
    Denote the quadratic form $q$ on $M \times_N M$ by the same letter ($q(m_1, m_2) = (q(m_1), q(m_2))$).

%
%

    \subsection{Invariance of $N_0$}
    From now on, assume that $(M_0, q_0)$ is a projective module of odd rank $n$ with a semiregular quadratic form.

    For example, it may happen that $(M_0, q_0)$ is split: $M_0 = f_{-\nu}K \oplus \ldots \oplus f_{-1}K \oplus f_{0}K \oplus f_1 K \oplus \ldots \oplus f_{\nu}K$, where $n = 2\nu + 1$, and $q_0(x_{-\nu}, \ldots, x_{-1}, x_0, x_1, \ldots, x_{\nu}) = x_0^2 + x_{-1}x_1 + \ldots + x_{-\nu}x_{\nu}$.

    In~\cite[1.2]{abe}, Abe proposed using the notion of an admissible pair — a pair of ideals $(J, I)$ in the ring $K$ satisfying the inclusions
    \[\all{I \subset J\\ 2J + J^{[2]} \subset I},\label{abe}\tag{$A$}\]
    where $J^{[2]} = \langle\defset{x^2}{x \in J}\rangle$.

    In the split case, there is a one-to-one correspondence between admissible pairs~\eqref{abe} and~\eqref{admissible}: to a pair of ideals $(J, I)$ associate the submodule $N_0 \coloneqq f_{-\nu}I \oplus \ldots f_{-1}I \oplus f_0 J \oplus f_1 I \oplus \ldots \oplus f_{\nu}I$ and the same ideal $I$.
    It turns out that the nonsplit case has no fundamental differences, as the following proposition states.
    \proposal{\label{admissible-pairs-correspondence}
    Let $(M_0, q_0)$ be a projective module of odd rank $n = 2\nu + 1$ with a semiregular quadratic form.

    Let $(N_0, I)$ be an admissible pair of the form~\eqref{admissible}.
    Set $J \coloneqq M_0^\lor N_0$.
    Then $(J, I)$ is an Abe admissible pair (of the form~\eqref{abe}), and this correspondence between pairs~\eqref{abe} and pairs~\eqref{admissible} is bijective.
    \prove{
        Construct the inverse map: to a pair~\eqref{abe} $(J, I)$ associate the pair consisting of the submodule $N_0 \coloneqq \defset{u \in M_0 J}{b_{q_0}(u, M_0) \subset I}$ and the ideal $I$.

        The constructed maps and the conditions on $(M_0, q_0)$ are preserved under localizations.
        The conditions~\eqref{abe} and~\eqref{admissible} are local.
        Therefore, to check that these maps are well-defined and mutually inverse, we can assume the ring $K$ is local.

        By Proposition~\ref{orthogonal-decomposition}, there is a decomposition $(M_0, q_0) = \bigboxplus\limits_{i = 1}^{\nu}(M_{i}', q_{i}') \boxplus (M_{0}', q_{0}')$, with $M_0' = \angles{f_0}$.

        Fix an admissible pair either of type~\eqref{abe} or of type~\eqref{admissible} and construct the other.
        In both cases it is easily verified that $N_0 = M_0 I + f_0 J$, which ensures that the constructed pair is admissible.
        This equality also shows that the maps are indeed mutually inverse.
%
%
%
    }
    }
    \corollary{\label{semiregular-then-invariant}
    Under these assumptions, $\O(M_0, q_0)N_0 \subset N_0$.
    \prove{
        Let $\phi \in \O(M_0, q_0)$.
        The pair $(\phi(N_0), I)$ is also admissible~\eqref{admissible}, and corresponds to the same pair~\eqref{abe} $(J, I)$ as the pair $(N_0, I)$.
        Thus, $\phi(N_0) = N_0$.
    }
    }

    \corollary{\label{factoring-through}
        $\O(M, q)N \le N$ (where, as before, $M = e_{-2}K \oplus e_{-1}K \oplus M_0 \oplus e_1 K \oplus e_2 K$, $(N_0, I)$ is an admissible pair, and $N =  e_{-2}I \oplus e_{-1}I \oplus N_0 \oplus e_1 I \oplus e_2 I$).
        \prove{
            The condition is local, so we can assume $K$ is a local ring.

            Let $\phi \in \O(M, q)$. Over a local ring, $\phi$ has a Gauss decomposition, and it suffices to check that $(\EO(M, q)\cdot D(M, q))N \le N$.
            This follows from conditions~\eqref{admissible} and Corollary~\ref{semiregular-then-invariant}.
        }
    }

    \subsection{Decomposing $\O(M \times_N M, q)$ into a Semidirect Product}
    Let $\phi \in \O(M, q)$.
    By Corollary~\ref{factoring-through}, if $(u_1, u_2) \in M \times_N M$, then $(\phi(u_1), \phi(u_2)) \in M \times_N M$ as well.
    In other words, the homomorphism $\Delta^{\O}: \O(M, q) \to \O(M \times_N M, q)$, $\phi \mapsto (\phi\pi_1, \phi\pi_2)$ is well-defined.

    This homomorphism has a retraction $\pi_1^{\O}: \O(M \times_N M, q) \to \O(M, q)$, $\phi \mapsto \pi_1 \phi \Delta$.
    Depict this in a diagram:
    \[\begin{tikzcd}
    {\O(M, q)}
          & {\O(M\times_{N}M, q)} & {\O(M, q)}
          \arrow["{\Delta^{\O}}", from=1-1, to=1-2]
          \arrow["\id", curve={height=-24pt}, from=1-1, to=1-3]
          \arrow["{\pi_1^{\O}}", from=1-2, to=1-3]
    \end{tikzcd}\]
    \lemma{
        For all $y \in K, y' \in K\times_{I}K$ and $u \in M_0, u' \in M_0 \times_{N_0} M_0$:
        \begin{align*}
            \Delta^{\O}(T_{ij}(y)) &= T_{ij}(\Delta(y)) \\
            \Delta^{\O}(T_{i}(u)) &= T_i(\Delta(u)) \\
            \pi_{1}^{\O}(T_{ij}(y')) &= T_{ij}(\pi_1(y')) \\
            \pi_{1}^{\O}(T_{i}(u')) &= T_i(\pi_1(u'))
        \end{align*}
    }
    \begin{corollary_env}
        The sequence of $\Delta^{\O}$ and $\pi_1^{\O}$ restricts to $\EO$:
        \[\begin{tikzcd}
        {\EO(M, q)}
              & {\EO(M\times_{N}M, q)} & {\EO(M, q)}
              \arrow["{\Delta^{\EO}}", from=1-1, to=1-2]
              \arrow["\id", curve={height=-24pt}, from=1-1, to=1-3]
              \arrow["{\pi_{1}^{\EO}}", from=1-2, to=1-3]
        \end{tikzcd}\]
    \end{corollary_env}

    Define the relative orthogonal subgroup in $\O(M, q)$
    \[\O(M, q, N_0, I) \coloneqq \defset{\phi \in \O(M, q)}{\forall u \in M: \phi(u) \equiv u \pmod{N}}\]
    In the case $N_0 = M_0 I$, this subgroup, of course, coincides with $\O(M, q, I)$.
    \lemma{
        \label{some-lemma}
        The embedding $\iota_2: \O(M, q, N_0, I) \to \Ker(\pi_{1}^{\O})$, $\phi \mapsto (\pi_1, \phi\pi_2)$ is an isomorphism.
        \prove{
            There is an inverse homomorphism $\iota_2^{-1}: \Ker(\pi_{1}^{\O}) \to \O(M, q, N_0, I)$, $\psi \mapsto \pi_2 \psi \Delta$.
        }
    }
    Define the relative elementary subgroup in $\O(M, q)$ \[\EO(M, q, N_0, I) \coloneqq \ind{}{\angles{T_{ij}(y), T_i(u)\mid y \in I, u \in N_0}}{^{\EO(M, q)}}.\]
    In the case $N_0 = M_0 I$, this subgroup, of course, coincides with $\ind{}{\mathrm{FO}(M, q, I)}{^{\EO(M, q)}}$.
%
    \lemma{\label{for-relativization}
    The image of $\EO(M, q, N_0, I)$ under $\iota_2$ is the kernel of $\pi_1^{\EO}$.
    \prove{
        Clearly, $\iota_2(\EO(M, q, N_0, I)) \le \Ker(\pi_{1}^{\EO})$.
        Check surjectivity.

        Represent elements of $\EO(M \times_{N} M, q)$ as pairs, embedding it into $\EO(M, q) \times \EO(M, q)$.
        Let $(\id, t) \in \Ker(\pi_{1}^{\EO})$.
        Since it belongs to $\EO$, it decomposes to a product $(\id, t) = {(s_1, t_1)}{} \proddots {(s_N, t_N)}$, where $s_i, t_i \in \EO(M, q)$ are elementary transvections of the same root type and $\forall m \in M: s_i(m) = t_i(m) \pmod{N}$. Hence it is clear that $t_i s_i^{-1} \in \EO(M, q, N_0, I)$ is some elementary transvection such that $\forall m \in M$: $t_i s_i^{-1}(m) \equiv m \pmod{N}$, i.e., a transvection of level $(N_0, I)$.
        Write
        \[(\id, t) = (\id, t_1 s_1^{-1}) (\id, s_1 t_2 s_2^{-1}s_1^{-1}) \proddots (\id, s_1 \proddots s_{N-1} t_N s_N^{-1} \proddots s_1^{-1}),\]
        which shows: indeed $t \in \EO(M, q, N_0, I)$.
    }
    }

    \subsection{Standard Commutator formulae}
    \lemma{\label{double-is-finitely-generated}
    Assume that $N_0$ is a finitely generated $K$-module.

    Then $M_0 \times_{N_0} M_0$ is a finitely generated $K \times_I K$-module, and the ideal $(q_0(w)\mid w \in M_0 \times_{N_0} M_0) = (1) \normeq K \times_I K$.
    \prove{
        There is an isomorphism of $K$-modules $M_0 \times_{N_0} M_0 \cong M_0 \oplus N_0$, so under the assumptions made, $M_0 \times_{N_0} M_0$ is a finitely generated $K$-module, hence certainly a finitely generated $K \times_I K$-module.

        By Lemma~\ref{unit-q}, there exist $w_1, \ldots, w_N \in M_0$ and $x_1, \ldots, x_N \in K$ such that $\sum\limits_{i=1}^{N}q_0(w_i)x_i = 1$.
        Under the diagonal embedding, $\sum\limits_{i=1}^{N}q_0(w_i, w_i)\cdot(x_i, x_i) = 1 \in K \times_I K$.
    }
    }
    \lemma{\label{construct-finitely-generated}
    Let $U \subset N_0$ and $Y \subset I$ be finite sets.
    Then there exist finitely generated $\tilde{N}_0 \le N_0$ and $\tilde{I} \subset I$ such that still $U \subset \tilde{N}_0$ and $Y \subset \tilde{I}$, and still $(\tilde{N}_0, \tilde{I})$ is an admissible pair.
    \prove{
        By enlarging $Y$, we can assume that $\forall u, u_1, u_2 \in U$: $q_0(u) \in Y$ and $b_{q_0}(u_1, u_2) \in Y$.
        Since $(N_0, I)$ is an admissible pair, still $Y \subset I$.

        Set $\tilde{I} \coloneqq \angles{Y}_K \normeq K$ and $\tilde{N}_0 \coloneqq \angles{U}_K + M_0 \tilde{I} \le M_0$.
        It is easy to check that now $(\tilde{N}_0, \tilde{I})$ is an admissible pair~\eqref{admissible} of finitely generated submodule and ideal.
    }
    }
    An admissible pair~\eqref{admissible} $(N_0, I)$, where $N_0$ is a finitely generated submodule and $I$ is a finitely generated ideal, will be called \emph{finitely generated}.
    \corollary{\label{finitely-generated-1}
    Let $\phi \in \EO(M, q, N_0, I)$.
    Then there exists a finitely generated admissible pair $(\tilde{N}_0, \tilde{I})$ such that $\phi \in \EO(M, q, \tilde{N}_0, \tilde{I})$.
    \prove{
        Write \[\phi = \ind{^{\psi_1}}{T_{ij}(y_1)}{} \cdot \ind{^{\psi_1'}}{T_{i}(u_1)}{} \proddots \ind{^{\psi_N}}{T_{ij}(y_N)}{} \cdot \ind{^{\psi_N'}}{T_{i}(u_N)}{}\]
        where $y_i \in I$, $u_i \in N_0$, $\psi_i, \psi_i' \in \EO(M, q)$.
        Apply Lemma~\ref{construct-finitely-generated} to $U = \{u_1, \ldots, u_N\}$ and $Y \coloneqq \{y_1, \ldots, y_N\}$.
    }
    }
    \corollary{\label{finitely-generated-2}
    Let $\phi \in \O(M, q, N_0, I)$.
    Then there exists a finitely generated admissible pair $(\tilde{N}_0, \tilde{I})$ such that $\phi \in \O(M, q, \tilde{N}_0, \tilde{I})$.
    \prove{
        Let $M_0 = \angles{u_1, \ldots, u_m}_K$.
        Let $V = \{e_{-2}, e_{-1}, u_1, \ldots, u_m, e_{1}, e_2\}$.
        In particular, $M = \angles{V}_K$.

        Let $w \in V$.
        It has a decomposition $w = e_{-2}\pi_{-2}(w) + e_{-1}\pi_{-1}(w) + \pi_0(w) + e_1 \pi_1(w) + e_2 \pi_2(w)$ with $\pi_i(w) \in I$ for $i \in \{\pm 1, \pm 2\}$ and $\pi_0(w) \in N_0$.
        Apply Lemma~\ref{construct-finitely-generated} to $U \coloneqq \defset{\pi_0(w)}{w \in V}$ and $Y = \defset{\pi_i(w)}{i \in \{\pm 1, \pm 2\}, w \in V}$.
    }
    }
    \proposal{\label{standard-commutator}
    Let $M_0$ be a projective module of odd rank with a semiregular quadratic form.
    For an admissible pair~\eqref{admissible} $(N_0, I)$:
    \bullets{
        \item $\EO(M, q, N_0, I) \normeq \O(M, q)$;
        \item $[\EO(M, q), \O(M, q, N_0, I)] \le \EO(M, q, N_0, I)$.
    }
    \prove{
        We will check on elements.
        By Lemmas~\ref{finitely-generated-1} and~\ref{finitely-generated-2}, we can reduce the pair $(N_0, I)$ to a finitely generated one.
        By Lemma~\ref{double-is-finitely-generated}, Theorem~\ref{eo-normality} is applicable, so\[\EO(M \times_{N} M, q) \normeq \O(M \times_{N} M, q)\tag{$\circ$}\label{double-normality}\]
        \bullets{
            \item Let $t \in \EO(M, q, N_0, I)$, $\phi \in \O(M, q)$.
            Conjugating $\iota_2(t)$ by the element $\Delta^{\O}(\phi)$, we get an element in $\iota_2(\O(M, q, N_0, I))$, since its first component is $\id$. By normality~\eqref{double-normality}, $\ind{^{\Delta^{\O}(\phi)}}{\iota_2(t)}{}$ also lies in $\EO(M \times_{N} M, q)$.
            By Lemma~\ref{for-relativization}, $\ind{^{\Delta^{\O}(\phi)}}{\iota_2(t)}{}$ lies in $\iota_2(\EO(M, q, N_0, I))$.
            Since $\O(M, q, N_0, I) \overset{\iota_2}\Map \O(M \times_{N} M, q)$ is injective, $\ind{^{\Delta^{\O}(\phi)}}{\iota_2(t)}{} = (\id, \ind{^{\phi}}{t}{})$ implies $\ind{^{\phi}}{t}{} \in \EO(M, q, N_0, I)$.
            \item Let $s \in \EO(M, q)$ and $t \in \O(M, q, N_0, I)$. In $ \O(M \times_{N} M, q)$:
            \[(s, s)(\id, t)(s^{-1}, s^{-1})(\id, t^{-1}) = (\id, [s, t]).\]
            It is clear that this lies in $\iota_2(\O(M, q, N_0, I)) = \Ker(\pi_{1}^{\O})$, and by normality~\eqref{double-normality} it also lies in $\EO(M \times_{N} M, q)$, so we get $[s, t] \in \EO(M, q, N_0, I)$.\qedhere
        }
    }
    }
    \section{Preliminary Computations}
    In this section, we assume that $M_0$ is a projective module of odd rank $n$ with a semiregular quadratic form, $(N_0, I)$ is an admissible pair, and $N = e_{-2}I \oplus e_{-1}I \oplus N_0 \oplus e_1 I \oplus e_2 I$.
    \subsection{$\EO$ is Perfect and the Relative Analogue}
    \proposal{\label{eo-perfection}
    Let $M_0$ be a projective module of odd rank $n \ge 3$ with a semiregular quadratic form, and $(N_0, I)$ an admissible pair.
    Then $\EO(M, q)$ is perfect.

    Moreover, $[\EO(M, q), \EO(M, q, N_0, I)] = \EO(M, q, N_0, I)$.
    \prove{
        Obviously follows from the commutation relations.
        Let's prove the relative result immediately. Denote $H \coloneqq [\EO(M, q), \EO(M, q, N_0, I)]$.

        Let $y \in I$. Since the pair $(N_0, I)$ is admissible, $M_0 I \le N_0$.
        For $i \ne \pm j$: $[T_{i}(u), T_{j}(vy)] = T_{-j,i}(b_{q_0}(u, v)y)$, and since by Lemma~\ref{unit-bq}, $(b_{q_0}(u, v)\mid u, v \in M_0) = (1) \normeq K$, we have $T_{-j,i}(y) \in H$.

        Further, $\forall u \in M_0$: $T_j(u) = [T_i(u), T_{ij}(1)]\cdot T_{-j,i}(-q_0(u)y) \in H$.
    }
    }
    \subsection{Computing Transporters}
    We compute $\Trans_{\O(M, q)}(\EO(M, q), \O(M, q, N_0, I)) = \defset{\phi \in \O(M, q)}{[\phi, \EO(M, q)] \le \O(M, q, N_0, I)}$.

    Clearly, $\O(M, q, N_0, I) \le \Trans_{\O(M, q)}(\EO(M, q), \O(M, q, N_0, I))$.
    \lemma{\label{factoring-module}
    The $K$-module $M/N$ has the structure of a $K/I$-module, and the form $q$ induces a form on $M/N$, and still $q_{M/N}$ has Witt index at least $2$.

    Moreover, $\rho_N: \O(M, q)/\O(M, q, N_0, I) \to \O(M/N, q_{M/N})$ is an embedding.
    \prove{
        Since the pair $(N_0, I)$ is admissible, $M/N$ is a $K/I$-module, and $q: M \to K$ factors through $N$, yielding $q_{M/N}: M/N \to K/I$.
        The direct summand $K^{\oplus 4}$ factors to a direct summand $(K/I)^{\oplus 4}$, on which the form is split.

        Since $(N_0, I)$ is admissible, by Corollary~\ref{factoring-through}, every $\phi \in \O(M, q)$ factors through the quotient, yielding $\rho_N(\phi) \in \O(M/N, q_{M/N})$.
        Finally, $\Ker(\rho_N) = \O(M, q, N_0, I)$ by the definition of the latter.
    }
    }
    Factoring the transporter by $\O(M, q, N_0, I)$, i.e., applying $\rho_N: \O(M, q) \to \O(M/N, q_{M/N})$, we get \[\Trans_{\O(M, q)}(\EO(M, q), \O(M, q, N_0, I)) = \rho_N^{-1}(\defset{\phi \in \O(M/N, q_{M/N})}{[\phi, \EO(M/N, q_{M/N})] = 1})\]
    \lemma{\label{commute}
    Assume that there are no conditions of projectivity $M_0$ or semiregularity of $q_0$.
    Even then, $\defset{\phi \in \O(M, q)}{[\phi, \EO(M, q)] = 1} = \Cent(\O(M, q)) = \mu_2(K)\cdot\id$.
    \prove{
        The inclusions $\ge$ are obvious, we prove \[\defset{\phi \in \O(M, q)}{[\phi, \EO(M, q)] = 1} \le \mu_2(K) \cdot \id.\]
        Let $[\phi, \EO(M, q)] = 1$.
        Since $[\phi, T_{ij}(1)] = 1$, and $[\phi, \id] = 1$, in $\End(M)$ we have the equality
        \[\phi(e_{ij} - e_{-j,-i}) = (e_{ij} - e_{-j,-i})\phi.\]
        For example, for $i = 1$, $j = 2$ the equality becomes
        \[\left(\arr{c|c|c|c|c}{&-[\phi]_{-2,-2} &&&[\phi]_{-2,1} \\\hline &-[\phi]_{-1,-2} &&&[\phi]_{-1,1} \\\hline &-[\phi]_{0,-2} &&&[\phi]_{0,1}\\\hline &-[\phi]_{1,-2} &&&[\phi]_{1,1} \\\hline &-[\phi]_{2,-2} &&&[\phi]_{2,1}}\right) =
        \left(\arr{c|c|c|c|c}{-[\phi]_{-1,-2} &-[\phi]_{-1,-1}&-[\phi]_{-1,0}&-[\phi]_{-1,1}&-[\phi]_{-1,2} \\\hline &&&& \\\hline &&&& \\\hline [\phi]_{2,-2} &[\phi]_{2,-1}&[\phi]_{2,0}&[\phi]_{2,1} &[\phi]_{2,2} \\\hline &&&&}\right).\]
        From this we get the equalities
        \begin{align*}
        [\phi]
            _{-j,-j} &= [\phi]_{-i,-i},\\
            -[\phi]_{-j,+i} &= [\phi]_{-i,+j},\\
            -[\phi]_{+i,-j} &= [\phi]_{+j,-i},\\
            [\phi]_{+i,+i} &= [\phi]_{+j,+j},\\
            \text{and all other }[\phi]_{-i,*} &= [\phi]_{j,*} = [\phi]_{*,-j} = [\phi]_{*,i} = 0.
        \end{align*}
        Substituting various pairs $i, j$, we get that $[\phi]$ has the form
        \begin{center}
            \scalebox{0.85}{
                \(
                \underset{i = 1, j = 2}{\left(\arr{c|c|c|c|c}{r_1 & * & * &y_1& * \\\hline & r_1 &&&-y_1 \\\hline &*&*&&* \\\hline y_2&*&*& r_2 &* \\\hline &-y_2&&& r_2}\right)};
                \underset{i = -1, j = -2}{\left(\arr{c|c|c|c|c}{r_3 &  &  &y_3&  \\\hline *& r_3 &*&*&-y_3 \\\hline *&&*&*& \\\hline y_4&&& r_4 & \\\hline *&-y_4&*&*& r_4}\right)};
                \underset{i = 1, j = -2}{\left(\arr{c|c|c|c|c}{r_5 &y_5 &&& \\\hline y_6 & r_6 &&& \\\hline *&*&*&& \\\hline *&*&*& r_5 &-y_5 \\\hline *&*&*&-y_6& r_6}\right)};
                \underset{i = -1, j = 2}{\left(\arr{c|c|c|c|c}{r_7 &y_7 &*&*&* \\\hline y_8 & r_8 &*&*&* \\\hline &&*&*&* \\\hline &&& r_7 &-y_7 \\\hline &&&-y_8& r_8}\right)}.
                \)}
        \end{center}
        Thus, \[[\phi] = \vect{r &&&& \\ &r &&& \\ &&\phi_0&& \\ &&&r& \\ &&&&r}.\]
        Since $\phi$ preserves the quadratic form, $r^2 = 1$.

        Since scalars lie in the center of $\O(M, q)$, $r^{-1}\phi$ also commutes with $\EO(M, q)$.
        Finally, $\ind{^{d(r^{-1}\phi_0)}}{T_{i}(u)}{} = T_i(r^{-1}\phi_0(u))$, and hence $r^{-1}\phi_0 = \id$, so $\phi = r\id$.
    }
    }
    \begin{proposal_env}
        \label{transporters}
        Let $M_0$ be a projective module of odd rank $n$ with a semiregular quadratic form, and $(N_0, I)$ an admissible pair.
        Then
        \begin{align*}
            \rho_N^{-1}(\mu_2(K/I)\cdot \id) =& \Trans_{\O(M, q)}(\EO(M, q), \O(M, q, N_0, I)), \\
            =& \Trans_{\O(M, q)}(\O(M, q), \O(M, q, N_0, I)),\\
            =& \Trans_{\O(M, q)}(\EO(M, q), \EO(M, q, N_0, I)).\\
        \end{align*}
        If $N_0 = M_0 I$, then we can say more:
        \[\rho_N^{-1}(\mu_2(K/I)\cdot \id) = \mu_2(K) \times \mathrm{SO}(M, q, N_0, I).\]
        \prove{
            Applying~\cref{commute} to the quotient module $M/N$ and quotient ring $K/I$, we get the equality for the first transporter.
            The second transporter is a subgroup of the first, and the explicit presentation of the first makes it clear they are equal.

            Now we verify that first and third transporters are equal.
            Let $T_N \coloneqq \Trans_{\O(M, q)}(\EO(M, q), \O(M, q, N_0, I))$.
            The result follows immediately from~\cref{eo-perfection,standard-commutator} and the Hall-Witt identity:
            \[[T_N, [\EO(M, q), \EO(M, q)]] \le [[T_N, \EO(M, q)], \EO(M, q)] \le \EO(M, q, N_0, I)\]

            If $N_0 = M_0 I$, then $N = MI$, and the quotient module $M/MI$ is a projective $K/I$ module.
            In this case, $\rho_N$ restricts to a homomorphism $\mathrm{SO}(M, q) \to \mathrm{SO}(M/MI, q_{K/I})$.
            So if $\phi\psi \in \rho_N^{-1}(\mu_2(K/I)\id_{M/N})$, where $\phi \in \mu_2(K)\id_M$ and $\psi \in \mathrm{SO}(M, q)$, then $\psi \in \Ker(\rho_N)$.
        }
    \end{proposal_env}
    There is also the transporter $\Trans_{\O(M, q)}(\O(M, q), \EO(M, q, N_0, I))$, and its explicit presentation is unclear.
    But we won't need it.

    \section{Proof of the Standard Classification}
    In this section, we assume that $M_0$ is a projective module of odd rank $n \ge 3$ with a semiregular quadratic form, $(N_0, I)$ is an admissible pair, and $N = e_{-2}I \oplus e_{-1}I \oplus N_0 \oplus e_1 I \oplus e_2 I$.

    \subsection{Constructing an Admissible Pair}
    \lemma{\label{admissible-pair-detection}
    Let $H \le \O(M, q)$ be an $\EO$-normal subgroup.
    For $i \in \{\pm 1, \pm 2\}$ define $N_0^{(i)} \coloneqq \defset{u \in M_0}{T_i(u) \in H}$ and for $i, j \in \{\pm 1, \pm 2\}$, $i \ne \pm j$ define $I^{(ij)} \coloneqq \defset{y \in K}{T_{ij}(y) \in H}$.

    Then $N_0 \coloneqq N_0^{(i)}$ does not depend on the choice of $i$, $I \coloneqq I^{(ij)}$ does not depend on the choice of $i$ and $j$, and $(N_0, I)$ is an admissible pair.
    \prove{
        Clearly, $N_0^{(i)} \subset M_0$ and $I^{(ij)} \subset K$ are additive subgroups.

        By Lemma~\ref{unit-bq}, there exist $u_k, v_k \in M_0, x_k \in K$ such that $\sum\limits_{k = 1}^{N}b_{q_0}(u_k, v_k)x_k = 1$. Let $y \in I^{(ij)}$ and $z \in K$. Then
        \multline{\prod\limits_{k = 1}^{N}[T_{i}(u_k + v_k x_k z), T_{ij}(y)] \cdot [T_{i}(-u_k), T_{ij}(y)] \cdot [T_{i}(-v_k x_k z), T_{ij}(y)] = \\=\prod\limits_{k = 1}^{N}T_{-j,i}(b_{q_0}(u_k, v_k x_k z)y) = T_{-j,i}(zy) \in H.}
        The calculation shows that $K \cdot I^{(ij)} \subset I^{(-j,i)} = I^{(-i,j)}$, and three other similar equalities imply that $I^{(ij)}$ is an ideal in $K$ independent of the choice of indices.

        Note that $P_{ij} \coloneqq T_{ij}(1)T_{ji}(-1)T_{ij}(1)$ has a monomial matrix. For example, for $i = 1, j = 2$:
        \[[P_{12}] = \smallvect {&-1&&&\\1&&&&\\&&1&&\\&&&&1\\&&&-1&}.\]

        Let $u \in N_0^{(i)}$, $v \in M_0$, $y \in I$, $z \in K$.
        \bullets{
            \item $P_{ij}^{-1} T_{i}(u) P_{ij} = T_{j}(u) \in H$, so $N_0^{(i)}$ is independent of the choice of index.
            \item $[T_i(u), T_{ij}(1)] = T_j(u) T_{-j,i}(q_0(u)) \in H$, so $q_0(N_0) \subset I$.
            \item$[T_{i}(u), T_{ij}(z)] = T_{j}(uz) T_{-j,i}(q_0(u)z) \in H$, so $N_0$ is a submodule.
            \item $[T_i(u), T_j(v)] = T_{-j,i}(b_{q_0}(u, v)) \in H$, so $b_{q_0}(N_0, M_0) \subset I$.
            \item $[T_{i}(v), T_{ij}(y)] = T_{j}(vy) T_{-j,i}(q_0(v)y) \in H$, so $M_0 I \le N_0$.\qedhere
        }
    }
    }

    \subsection{Extracting Transvections}
    \lemma{\label{commutator-identity}
    Let $G$ be a group, $x, y, z \in G$. Then $\ind{^x}{[y, x^{-1}z]}{} = [x, y]\cdot[y, z]$.
    }
    \lemma{\label{commutator-distributivity}
    Let $G$ be a group, $x, y, z \in G$. Then $[xy, z] = \ind{^x}{[y, z]}{} \cdot [x, z]$.
    In particular, if $[y, z] = 1$, then $[xy, z] = [x, z]$.
    }
    \lemma{\label{bottom-unipotent}
    Let $\phi \in \O(M, q)$ have a matrix of the form
        $[\phi] = \smallvect{1&&&& \\ *&1&&& \\ *&&1&& \\ *&&&1& \\ *&*&*&*&1}$.

        Then $\phi = T_{-1,-2}([\phi]_{-1,-2}) T_{-2}([\phi]_{0,-2}) T_{1,-2}([\phi]_{1,2})$.
        \prove{
            Set $\psi \coloneqq T_{-1,-2}([\phi]_{-1,-2}) T_{-2}([\phi]_{0,-2}) T_{1,-2}([\phi]_{1,2})$.

            Then
            $[\psi^{-1}\phi] = \smallvect{1&&&& \\ &1&&& \\ &&1&& \\ &&&1& \\ *&*&*&*&1}$,
            and by Lemma~\ref{bottom-row-is-zero}, we get $\psi^{-1}\phi = \id$.
        }
    }
    From now on, assume that $H \le \O(M, q)$ is an $\EO$-normal subgroup, and $\EO(M, q, N_0, I) \le H$.

    The proof of the following lemma imitates computations in~\cite[5.22.i]{preusser}.
    Unfortunately, in the nonsplit case these computations have become somewhat longer.
    \lemma{\label{calculation}
    Let $\phi \in H$ be such that $[\phi]_{1,2} \notin I$.
    Then $H$ contains an elementary transvection $T_{ij}(z')$, where $z' \notin I$.
    \prove{Let $u \in M_0$ be some vector, to be chosen later.
    Set \gather{\tau \coloneqq T_2(u[\phi]_{1,2})T_1(-u[\phi]_{2,2})  T_{-1,2}(-u^*[\phi]_{0,2} + [\phi]_{1,2}[\phi]_{2,2}q_0(u))\\
    [\tau] = \left(\arr{c|c|c|c|c} {1&&-[\phi]_{1,2}u^*&u^*[\phi]_{0,2} + [\phi]_{1,2}[\phi]_{2,2}q_0(u)&-[\phi]_{1,2}^2 q_0(u) \\\hline {}&1&[\phi]_{2,2}u^*&-[\phi]_{2,2}^2 q_0(u)&-u^*[\phi]_{0,2} + [\phi]_{1,2}[\phi]_{2,2}q_0(u) \\\hline {}&&1&-u[\phi]_{2,2}&u[\phi]_{1,2} \\\hline {}&&&1& \\\hline {}&&&&1}\right)}
    A direct calculation shows that $[\tau \phi]_{*,2} = [\phi]_{*,2}$. Set $\xi \coloneqq \phi^{-1}\tau \phi$, clearly $[\xi]_{*,2} = e_{2}$. By Lemma~\ref{bottom-row-is-zero} for $e_2$ instead of $e_{-2}$,
        $[\xi] = \smallvect{1&&&&\\ *&*&*&*& \\ *&*&*&*& \\ *&*&*&*& \\ *&*&*&*&1}$.
        By the same Lemma~\ref{bottom-row-is-zero}, the matrix $\xi^{-1}$ also has this form.

        Now set
        \[\zeta \coloneqq \ind{^{\tau}}{[T_{2,1}(1), [\tau^{-1}, \phi^{-1}]]}{} = \ind{^\tau}{[T_{2,1}(1), \tau^{-1}\xi]}{} \underset{\ref{commutator-identity}}= [\tau, T_{2,1}(1)] \cdot [T_{2,1}(1), \xi]\]
        Denote $z \coloneqq [\phi]_{1,2}$, and compute the first commutator factor:
        \[[\tau, T_{2,1}(1)] \underset{\ref{commutator-distributivity}}= [T_2(u z), T_{2,1}(1)] = T_{1}(u z)T_{-1,2}(q_0(u)z^2).\]
        To compute the second commutator factor, note that
        \[\left[\ind{^{T_{2,1}(1)}}{\xi}{}\right] = \left(\arr{c|ccc|c}{1&&&& \\\hline {}*&[\xi]_{-1,-1}&[\xi]_{-1,0}&[\xi]_{-1,1}& \\ {}*&[\xi]_{0,-1}&[\xi]_{0,0}&[\xi]_{0,1}& \\ *&[\xi]_{1,-1}&[\xi]_{1,0}&[\xi]_{1,1}& \\\hline * & * & * & * &1}\right),\]
        so the commutator $[T_{2,1}(1), \xi]$ has the form $\smallvect{1&&&& \\ *&1&&& \\ *&&1&& \\ *&&&1& \\ *&*&*&*&1}$; by Lemma~\ref{bottom-unipotent}, $[T_{2,1}(1), \xi] = T_{-1,-2}(*)T_{-2}(*)T_{1,-2}(*)$.
        Let $v \in M_0$ be another vector, also to be chosen later. Commute \multline{\eta \coloneqq [\zeta, T_{-2}(v)] \underset{\ref{commutator-distributivity}}= [T_{1}(uz) T_{-1,2}(q_0(u)z^2), T_{-2}(v)] \underset{\ref{commutator-distributivity}}=\\= \ind{^{T_1(uz)}}{[T_{-1,2}(q_0(u)z^2), T_{-2}(v)]}{} \cdot [T_{1}(uz), T_{-2}(v)] =\\= T_{2,1}(-q_0(u)q_0(v)z^2)T_{1}(vq_0(u)z^2) \cdot T_{2,1}(b_{q_0}(u, v)z).}
        \bullets{
            \item If $z^2 \in I$, then it suffices to choose $u, v \in M_0$ such that $b_{q_0}(u, v)z \notin I$ (which is possible by Lemma~\ref{unit-bq}).
            Then $\eta \equiv T_{2,1}(b_{q_0}(u, v)z) \pmod{\EO(M, q, N_0, I)}$.
            \item If $z^2 \notin I$, then we need to choose $u, v, w \in M_0$ such that $b_{q_0}(w, v) q_0(u) z^2 \notin I$ (which is possible by Lemmas~\ref{unit-q} and~\ref{unit-bq}).
            Then \[[T_{-2}(w), \eta] = T_{-1,-2}(b_{q_0}(w, v)q_0(u)z^2)\qedhere\]
        }
    }
    }
    \lemma{\label{common-entries-in-ideal}
    Let $\phi \in H$ and for some distinct $i, j \in \{\pm 1, \pm 2\}$: $[\phi]_{ij} \notin I$.
    Then $H$ contains an elementary transvection of the form $T_{i'j'}(z)$, where $z \notin I$.
    \prove{
        If $i \ne -j$, then this is Lemma~\ref{calculation} up to index change.

        If $i = -j$, then we can choose $k \in \{\pm 1, \pm 2\}$ such that $k \ne \pm i$.
        Then $\ind{^{T_{k,i}(1)}}{\phi}{}$ at position $(k, -i)$ has $[\phi]_{i,-i} + [\phi]_{k,-i}$.
        Clearly, $\any{[\phi]_{i,-i}  + [\phi]_{k,-i}\notin I \\ [\phi]_{k,-i} \notin I}$, and depending on which case holds, apply Lemma~\ref{calculation} to $\ind{^{T_{k,i}(1)}}{\phi}{}$ or to $\phi$.
    }
    }
%
    Depending on context, an overline will denote the image of a ring element $K$ in the quotient ring $K/I$, of a module element $M_0$ in the quotient module $M_0/N_0$, or the image of an element of $\O(M_0, q_0)$ under the (not necessarily surjective) homomorphism $\O(M_0, q_0) \to \O(M_0/N_0, (q_0)_{M_0/N_0})$
    \lemma{\label{common-entries-in-submodule}
    Let for every $\phi \in H$, for all distinct $i, j \in \{\pm 1, \pm 2\}$: $[\phi]_{ij} \in I$.
    Then the image of $\phi$ in $\O(M/N, q_{M/N})$, call it $\overline{\phi}$, has the form
        \[[\overline{\phi}] = \vect{\overline{r} &&&& \\ &\overline{r}&&& \\ \overline{v}_{-2} & \overline{v}_{-1} & \overline{\phi}_0 & \overline{v}_1 & \overline{v}_2 \\ &&&\overline{r}& \\ &&&&\overline{r}},\]
        where $\overline{r} \in \mu_2(K/I)$, $\overline{v}_j \in M_0/N_0$, with $\overline{v}_j^* = 0$ and $(q_0)_{M/N}(\overline{v}_j) = 0$.
        \prove{
            It suffices to assume $I = 0$ and $N_0 = 0$, since we can homomorphically map $H$ to $\O(M/N, q_{M/N})$, and it will remain $\EO$-normal.

            By assumption, ${\phi} \in \O(M, q)$ has the form
            \[[\phi] = \vect{r_{-2} &&{\theta}_{-2}&& \\ &r_{-1}&{\theta}_{-1}&& \\ {v}_{-2} & {v}_{-1} & {\phi}_0 & {v}_1 & {v}_2 \\ &&{\theta}_1& r_1 & \\ &&{\theta}_2&& r_2},\]
            where $r_i \in K$, ${v}_j \in M_0$, ${\theta}_i \in M_0^\lor$.

            Let ${w} \in M_0$ be an arbitrary vector, and $i, j \in \{\pm 1, \pm 2\}$ such that $i \ne \pm j$.
            \bullets{
                \item Note that $\left[\ind{^{T_{i}({w})}}{{\phi}}{}\right]_{-i,j} = -{w}^* {v}_j$.
                But by assumption this is $0$, and hence all ${v}_j^* = 0 \in M_0^\lor$.

                \item Note that $\left[\ind{^{T_{ij}(1)}}{\phi}{}\right]_{ij} = r_{j} - r_{i}$, so $r_i = r_j$.
                Substituting $-i$ for $i$, we get that all $r_i$ are equal to each other and equal to some $r \in K$.
                Since $\phi$ preserves the quadratic form, $0 = q(\phi(e_i)) = q_0(v_i)$, so $r^2 = 1$.
                \item Generally, $\phi({w}) = e_{-2}{\theta}_{-2}(w) + e_{-1}{\theta}_{-1}(w) + {\phi}_0(w) + e_1 {\theta}_1(w) + e_2 {\theta}_2(w)$.
                But ${\phi}(w)$ is orthogonal to all $\phi(e_i)$, hence $0 = r_i\theta_{-i}(w) + \phi_0(v_{i})^*(\phi_0(w)) = r_i\theta_{-i}(w)$, hence $\theta_{-i} = 0$.\qedhere
            }
        }
    }
    \lemma{\label{common-entries-in-submodule-2}
    Let for every $\phi \in H$, for all distinct $i, j \in \{\pm 1, \pm 2\}$: $[\phi]_{ij} \in I$.
    If there exists $\phi \in H$ and $j \in \{\pm 1, \pm 2\}$ such that $[\phi]_{0,j} \notin N_0$, then there is also a transvection of the form $T_{i}(v) \in H$, where $v \notin N_0$.
    \prove{
        By Lemma~\ref{common-entries-in-submodule}, all $[\phi]_{0,i}$ are such that $b_{q_0}([\phi]_{0,i}, M_0) \in I$ and $q_0([\phi]_{0,i}) \in I$.
        So we can write for some $\overline{v}_j$, $\overline{r}$, $\overline{\phi}_0$ the equality in $\O(M/N, q_{M/N})$: \[\overline{\phi} =  \underbrace{T_{-2}(\overline{v}_{-2}\overline{r})T_{-1}(\overline{v}_{-1}\overline{r})T_{1}(\overline{v}_1 \overline{r})T_2(\overline{v}_2 \overline{r})}_{\tau}\cdot \underbrace{\smallvect{\overline{r}&&&& \\ &\overline{r} &&& \\ &&\overline{\phi}_0 && \\ &&& \overline{r} & \\ &&&& \overline{r}}}_{d}.\]
        From this representation, it is clear that \[[\overline{\phi}, T_{ji}(1)] \underset{\ref{commutator-distributivity}}=  [\tau, T_{ji}(1)] = T_i(\overline{v}_j \overline{r}) T_{-j}(\overline{v}_{-i}\overline{r}) \in \O(M/N, q_{M/N}).\]
        Lift this equality to $\O(M, q)$ as
        \[[\phi, T_{ji}(1)] \equiv T_{i}({v}_j r)T_{-j}({v}_{-i} r)\pmod{\O(M, q, N_0, I)}, \]
        and commute it with $T_{i,-j}(r)$:
        \[[[\phi, T_{ji}(1)], T_{i,-j}(r)]\equiv T_{-j}(v_j r^2) \equiv T_{-j}(v_j);\]
        here the equality is modulo $[\O(M, q, N_0, I), \EO(M, q)] = \EO(M, q, N_0, I)$.
        Since $\EO(M, q, N_0, I) \le H$, we have $T_{-j}(v_j) \in H$.
    }
    }

    \lemma{\label{dealing-with-center}
    Let for every $\phi \in H$, for all distinct $i, j \in \{\pm 1, \pm 2\}$: $[\phi]_{ij} \in I$ and $[\phi]_{0,j} \in N_0$.
    Then $\phi \in \Trans_{\O(M, q)}(\EO(M, q), \O(M, q, N_0, I))$.
    \prove{Map $H$ homomorphically to $\O(M/N, q_{M/N})$. Under this procedure, it will remain $\EO$-normal, and by~\ref{transporters}, we need to ensure that $H$ consists of scalar transformations.

    Let $\phi \in H$ and $w \in M_0$ be arbitrary.
    By Lemmas~\ref{common-entries-in-submodule} and~\ref{common-entries-in-submodule-2}, $\phi = \smallvect{r &&&& \\ & r &&& \\ && \phi_0 && \\ &&& r & \\ &&&& r}$.
    Compute $\left[\ind{^{T_{j}(w)}}{\phi}{}\right]_{0,j} = w r - \phi_0(w)$; by assumption this equals $0 \in M_0/N_0$.
    In other words, $\phi_0 = r\id$ in $\End(M_0/N_0)$.
    }
    }
    \theorem{\label{main-theorem}
    Let $(M_0, q_0)$ be a projective module of odd rank $n \ge 3$ with a semiregular quadratic form.
    Let $H \le \O(M, q)$ be a subgroup.

    Then $H$ is $\EO$-normal if and only if there exists an admissible pair~\eqref{admissible} $(N_0, I)$ such that \[\EO(M, q, N_0, I) \le H
    \le \Trans_{\O(M, q)}(\EO(M, q), \O(M, q, N_0, I)).\]
        Moreover, if $H$ is $\EO$-normal, then such a pair $(N_0, I)$ is unique.
        \prove{
            Denote $T_{N} \coloneqq \Trans_{\O(M, q)}(\EO(M, q), \O(M, q, N_0, I))$.

            Any $H \le \O(M, q)$ contained between $\EO(M, q, N_0, I)$ and $T_N$ is $\EO$-normal, since by~\cref{transporters} $[H, \EO(M, q)] \le \EO(M, q, N_0, I) \le H$.

            Set $N_0 \coloneqq \defset{u \in M_0}{\exists i \in \{\pm 1, \pm 2\}: T_i(u) \in H}$ and $I \coloneqq \defset{y \in K}{\exists i, j \in \{\pm 1, \pm 2\}: T_{ij}(y) \in H}$.
            By Lemma~\ref{admissible-pair-detection}, $(N_0, I)$ is an admissible pair, and $\EO(M, q, N_0, I) \le H$. As usual, denote $N \coloneqq e_{-2}I \oplus e_{-1}I \oplus N_0 \oplus e_1 I \oplus e_2 I$.

            Now we show that $H \le T_N$.
            By Lemma~\ref{common-entries-in-ideal}, $\forall \phi \in H$, for all distinct $i, j \in \{\pm 1, \pm 2\}$: $[\phi]_{ij} \in I$. So Lemma~\ref{common-entries-in-submodule-2} applies, and $[\phi]_{0,j} \in N_0$.
            Hence Lemma~\ref{dealing-with-center} applies, and $\phi \in T_N$.

            Uniqueness follows for the following reason: if $\EO(M, q, N_0, I) \le H \le T_N$, then all elementary transvections with parameters of level $(N_0, I)$ are contained in $H$, and the explicit form of the transporter~\ref{transporters} shows that elementary transvections with parameters $u \notin N_0$ or $y \notin I$ are not contained in $H$.
        }
    }
    \section*{Acknowledgments}
    I want to thank Egor Voronetsky for his scientific supervision and guidance throughout this work.
    
\end{document}